\documentclass{article}
\usepackage{amssymb}

\usepackage{amsmath}

\usepackage{amsthm}

\newtheorem{thm}{Theorem}[section]

\newtheorem{lem}[thm]{Lemma}

\newtheorem{rem}[thm]{Remark}
\newtheorem{quest}[thm]{Question}

\title{Cubic Time Recognition of Cocircuit Graphs of Uniform Oriented Matroids}

\author{Stefan Felsner, Ricardo G{\'o}mez, Kolja Knauer, \\Juan Jos{\'e} Montellano-Ballesteros, Ricardo Strausz}

\begin{document}

\maketitle

\begin{abstract}
We present an algorithm which takes a graph as input and decides in cubic 
time if the graph is the cocircuit graph of a uniform oriented matroid. In the
affirmative case the algorithm returns the set of signed cocircuits of the 
oriented matroid. This improves an algorithm proposed by Babson, Finschi and Fukuda.

Moreover we strengthen a result of Montellano-Ballesteros and Strausz about
\textit{crabbed} connectivity of cocircuit graphs of uniform oriented matroids.
\end{abstract}

\section{Introduction}\label{intr}

The cocircuit graph is a natural combinatorial object associated with oriented matroids. In the case of spherical pseudoline-arrangements, i.e., rank $3$ oriented matroids, its vertices are the intersection points of the
lines and two points share an edge if they are adjacent on a line.
More generally, the Topological Representation Theorem of Folkman and Lawrence~\cite{Fol-78} says that every 
oriented matroid can be represented as an arrangement of pseudospheres. The cocircuit graph is the $1$-skeleton of this arrangement.

Cordovil, Fukuda and Guedes de Oliveira~\cite{Cor-00} show that 
a \textit{uniform} oriented matroid is determined by its cocircuit graph together with an \textit{antipodal labeling }.
Babson, Finschi and Fukuda~\cite{Bab-01} provide a polynomial time recognition algorithm for cocircuit graphs
of uniform oriented matroids, which reconstructs a uniform oriented matroid from its cocircuit graph up to isomorphism.

In~\cite{Mon-06}, Montellano-Ballesteros and Strausz provide 
a characterization of uniform oriented matroids in view of \textit{sign labeled} cocircuit graphs. We prove
a stronger version of that characterization in Theorem~\ref{thm:main}.

After introducing basic notions of oriented matroids we describe an algorithm which, 
given a graph $G$, decides in polynomial time if $G$ is the cocircuit graph of a uniform oriented matroid. In the
affirmative case the algorithm returns the set of signed cocircuits of the oriented matroid. Making use of Theorem~\ref{thm:main}
we obtain a better runtime than~\cite{Bab-01}. This in particular answers a question asked by Babson et al. However, some parts of 
our algorithm are identical to the one in~\cite{Bab-01}. At the 
end we will look at another question posed in~\cite{Bab-01} concerning
\textit{antipodality} in cocircuit graphs. We support the feeling that the antipodality problem is deep and hard
 by relating it to the Hirsch conjecture (cf.~\cite{Zie-95}).

\section{Properties of Cocircuit Graphs}

The notions introduced here are specialized to uniform oriented matroids, for a more general introduction see~\cite{Bjo-93}. 

A \textit{sign vector} on a ground set $E$ is an $X\in\{+,-,0\}^E$. The \textit{support of $X$} is $\underline{X}:=\{e\in E\mid X_e\neq 0\}$. By $X^0$ we denote the \textit{zero-support} $E\backslash\underline{X}$. The \textit{separator of two sign vectors $X,Y$} is defined as $S(X,Y):=\{e\in E\mid \{X_e,Y_e\}=\{+,-\}\}$. For a sign vector $X$ the sign vector $-X$ is the one where all signs are reversed.

We define a \textit{uniform oriented matroid} of \textit{rank} $r$ as a pair $\mathcal{M}=(E,\mathcal{C}^*)$ where $E$ is the ground set and the set $\mathcal{C}^*\subseteq\{+,-,0\}^E$ are the \textit{cocircuits} of $\mathcal{M}$. Denote by $n$ the size of $E$. Then $\mathcal{C}^*$ must satisfy the following axioms$\colon$
\begin{itemize}
 \item[(C1)] Every $X\in\mathcal{C}^*$ has support of size $n-r+1$.
\item[(C2)] For every $I\subseteq E$ of size $n-r+1$ there are exactly two cocircuits $X,Y$ with support $I$. Moreover $-X=Y$.
\item[(C3)] For every $X,Y\in\mathcal{C}^*$ and $e\in S(X,Y)$ there is a $Z\in\mathcal{C}^*$ with $Z_e=0$ and $Z_f\in\{X_f,Y_f,0\}$, for every $f\in E\backslash\{e\}$.
\end{itemize}

In the rest of the paper we will abbreviate uniform oriented matroid by {\rm{OM}}. On the set $\mathcal{C}^*$ of cocircuits of an {\rm{OM}} $\mathcal{M}$ one defines the \textit{cocircuit graph} $G_{\mathcal{M}}$ by making $X$ and $Y$ adjacent if they differ only a ``little bit'', i.e., $|X^0\Delta Y^0|=2$ and $S(X,Y)=\emptyset$. 

A more general notion is the following. Given a graph $G=(V,\mathcal{E})$ with vertices $V$ and edges $\mathcal{E}$ let $L:V\rightarrow \mathcal{S}$ be a bijection to a set of sign vectors $\mathcal{S}$ on a ground set $E$, which satisfies axioms (C1) and (C2). We call $L$ a \textit{sign labeling} of $G$ if we have $\{v,w\}\in \mathcal{E}$ if and only if $|L^0(v)\Delta L^0(w)|=2$ and $S(L(v),L(w))=\emptyset$. Every sign labeling $L$ of $G$ comes with the two parameters $r$ and $n$. 

If $G$ has a sign labeling with $\mathcal{S}$ satisfying also axiom (C3), i.e. $\mathcal{C}^*:=\mathcal{S}$ is the set of cocircuits of a {\rm{OM}} $\mathcal{M}$, then $G$ is a cocircuit graph $G_{\mathcal{M}}$. We then call $L$ an \textit{OM-labeling}. In~\cite{Bab-01} it is shown that $G_{\mathcal{M}}\cong G_{\mathcal{M}'}$ if and only if $\mathcal{M}\cong\mathcal{M}'$.

Clearly, $G_{\mathcal{M}}$ has exactly $2\binom{n}{n-r+1}$ vertices. In this section, such as the following, all the lemmas are well-known.
\begin{lem}\label{lem:reg}
Let $G_{\mathcal{M}}=(V,\mathcal{E})$ be a cocircuit graph with OM-labeling $L$ and $v\in V$. Then for every $f\in L^0(v)$ there are exactly two neighbors $u,w$ of $v$ with $L(u)_f=-$ and $L(w)_f=+$. In particular $G_{\mathcal{M}}$ is $2(r-1)$-regular.
\end{lem}
\begin{proof}
Let $L$ be an OM-labeling $v$ a vertex and $f\in L^0(v)$. Let $w$ be a vertex with $L(w)_f=+$ and $\underline{L(w)}\backslash\{f\}\subseteq\underline{L(v)}$, such that $S(L(v),L(w))$ is minimal. If there is an $e\in S(L(v),L(w))$ then we apply (C3) to $L(v), L(w)$ with respect to $e$. We obtain $L(u)$ with $\underline{L(u)}\backslash\{f\}\subseteq\underline{L(v)}\backslash\{e\}$ and $L(u)_f=+$. Since $L(u)_e=0$ and $L(u)_g\in\{L(v)_g,L(w)_g,0\}$ for $g\neq e$  we have $S(L(v),L(u))\subset S(L(v),L(w))$, a contradiction. Thus, $w$ is adjacent to $v$. If there was another neighbor $w'$ of $v$ with $L(w')_f=+$ then (C3) applied to $L(w')$ and $-L(w)$ with respect to $f$ would yield a cocircuit $L(u)$ with $\underline{L(u)}=\underline{L(v)}$. It is easy to see that $L(u)\neq\pm L(v)$, a contradiction to (C2).
\end{proof}

A contraction minor of a {\rm{OM}} $\mathcal{M}=(E,\mathcal{C}^*)$ is a {\rm{OM}} of the form $\mathcal{M}/E'=(E\backslash E',\mathcal{C}^*/E')$ where $E'\subseteq E$ and $\mathcal{C}^*/E':=\{X_{E\backslash E'}\mid X\in\mathcal{C}^*, E'\subseteq X^0\}$. Here $X_{E\backslash E'}$ denotes the restriction of $X$ to coordinates ${E\backslash E'}$. Generally $G_{\mathcal{M}/E'}$ is an induced subgraph of $G_{\mathcal{M}}$. The rank of $\mathcal{M}/E'$ is $r-|E'|$. If $\mathcal{M}/E'$ has rank $2$ we call it a \textit{coline} of $\mathcal{M}$.

\begin{lem}\label{lem:rk2}
Let $G_{\mathcal{M}}$ be a cocircuit graph with OM-labeling $L$ and $v,w\in V$ with $L(v)\neq -L(w)$. If $L(v)$ and $L(w)$ lie in a coline $\mathcal{M}'$ then $d(v,w)=|S(L(v),L(w))|+\frac{1}{2}|L^0(v)\Delta L^0(w)|$ and the unique $(v,w)$-path of this length lies in $G_{\mathcal{M}'}$.
\end{lem}
\begin{proof}
First note that $|S(L(v),L(w))|+\frac{1}{2}|L^0(v)\Delta L^0(w)|$ is a lower bound for the distance in any sign labeled graph. These are just the necessary changes to transform one sign vector into the other. To see that in a coline there exists a path of this length we use induction on $|S(L(v),L(w))|+\frac{1}{2}|L^0(v)\Delta L^0(w)|$. The induction base is clear per definition of sign labeling. So we proceed with the induction step.

Since $v,w$ lie in a coline we have $L^0(v)\backslash L^0(w)=\{e\}$ for some $e\in E$. By Lemma~\ref{lem:reg} vertex $v$ has a unique neighbor $u$ with $L(u)_e=L(w)_e$. If $S(L(u),L(w))=\emptyset$ we have $u=w$, because otherwise $v$ would have two neighbors with $L(u)_e=L(w)_e$, a contradiction to Lemma~\ref{lem:reg}. If $L(u)_f=0$ but $L(v)_f=L(w)_f\neq 0$ for some $f\in E$ then we can apply (C3) to $L(w)$ and $-L(u)$ with respect to $e$. Any resulting vector has support $\underline{L(v)}$. Since $S(L(u),L(w))\neq\emptyset$ it cannot have sign labeling $-L(v)$. On the other hand its $f$-entry equals $-L(v)_f$, i.e., it cannot have sign labeling $L(v)$ either, a contradiction to (C2). This yields $|S(L(u),L(w))|=|S(L(v),L(w))|-1$. Morover $|L^0(v)\Delta L^0(w)|$ is not decreased, since $v,w$ are not adjacent. Applying induction hypothesis gives the result. 

For a neighbor $u'$ of $v$ not in a coline with $v,w$ it is easy to check that $|S(L(u'),L(w))|+\frac{1}{2}|L^0(u')\Delta L^0(w)|\geq|S(L(v),L(w))|+\frac{1}{2}|L^0(v)\Delta L^0(w)|$. Hence $u'$ cannot lie on a shortest $(v,w)$-path.
\end{proof}

Let $G$ be a graph with sign labeling $L:V\rightarrow\mathcal{S}$. Let $X,Y\in \mathcal{S}$. We say that a path $P$ in $G$ is \textit{$(X,Y)$-crabbed} if for every vertex $w\in P$ we have $L(w)^+\subseteq X^+\cup Y^+$
and $L(w)^-\subseteq X^-\cup Y^-$. We call a $(u,v)$-path just \textit{crabbed} if it is $(L(u),L(v))$-crabbed. 
The following theorem is a strengthening of the main result of~\cite{Mon-06}$\colon$
\begin{thm}\label{thm:main}
 Let $L$ be a sign labeling of $G$. Then the following are equivalent:
\begin{enumerate}
 \item $L$ is an OM-labeling.
\item For all $v,w\in V$ there are exactly $|L(v)^0\backslash L(w)^0|$ vertex-disjoint crabbed $v,w$-paths.
\item For all $v,w\in V$ with $L(v)^0\neq L(w)^0$ there exists a crabbed $v,w$-path.
\end{enumerate}
\end{thm}
\begin{proof}
$(1)\Rightarrow (2):$ Let $G=G_{\mathcal{M}}$ a cocircuit graph. First by Lemma~\ref{lem:reg} it is clear that between any two vertices $v,w$ there can be at most $|L(v)^0\backslash L(w)^0|$ vertex-disjoint crabbed $v,w$-paths. For the other inequality we proceed by induction on the size of the ground set $E$ of $\mathcal{M}$. The induction base is skipped. For the inductive step we have to distinguish three cases.

If there is some $e\in L(v)^0\cap L(w)^0$ then consider the contraction minor $\mathcal{M}/\{e\}$. By induction hypothesis there are at least $|(L(v)^0\backslash\{e\})\backslash(L(w)^0\backslash\{e\})|=|L(v)^0\backslash L(w)^0|$ vertex-disjoint crabbed $v,w$-paths in $G_{\mathcal{M}/\{e\}}$. Since the latter is an induced subgraph of $G_{\mathcal{M}}$ we are done.

Otherwise, if $S(L(v),L(w))=\emptyset$ then $L(v),L(w)$ lie in a \textit{tope} of $\mathcal{M}$. This is the set of sign vectors $X$ with $X^+\subseteq L(v)^+\cup L(w)^+$ and $X^-\subseteq L(v)^-\cup L(w)^-$. Topes of a rank $r$ \rm{OM} are $(r-1)$-dimensional PL-spheres and hence their graph is $(r-1)$-connected~\cite{Bjo-93}. A less topological argument for the same fact can be found in~\cite{Cor-93}.

Otherwise, if there is some $e\in S(L(v),L(w))$ we consider the \textit{deletion minor} $\mathcal{M}\backslash\{e\}$. It is the oriented matroid on the ground set $E\backslash\{e\}$ with cocircuit set $\mathcal{C}^*\backslash\{e\}:=\{X_{E\backslash\{e\}}\mid X\in\mathcal{C}^*, X_{e}\neq 0\}$. By induction hypothesis there are $|L(v)^0\backslash L(w)^0|$ vertex-disjoint crabbed $v,w$-paths in $G_{\mathcal{M}\backslash\{e\}}$. If on such a path $P$ in $G_{\mathcal{M}\backslash\{e\}}$ two consecutive vertices $x,y$ have $e=S(L(x),L(y))$ then we apply (C3) with respect to $e$. We obtain a unique vertex $z$ with $L(z)_e=0$ and $L(z)_f=L(x)_f$ if $L(x)_f\neq 0$ and $L(z)_f=L(y)_f$ otherwise. The new vertex $z$ is adjacent to $x,y$ in $G_{\mathcal{M}}$. In this way we can extend $P$ to a crabbed $(v,w)$-path $P'$ in $G_{\mathcal{M}}$. 

Now suppose two different extended paths $P_1'$ and $P_2'$ share a vertex $z$. Thus there are mutually different $x_1,y_1,x_2,y_2$ yielding $z$. This implies that their labels $L(x_1)_{E\backslash\{e\}},L(y_1)_{E\backslash\{e\}},L(x_2)_{E\backslash\{e\}},L(y_2)_{E\backslash\{e\}}$ in $\mathcal{M}\backslash\{e\}$ have mutually empty separator. Moreover the zero-supports of these labels have mutually symmetrical difference two. Hence in $G_{\mathcal{M}\backslash\{e\}}$ the vertices $x_1,y_1,x_2,y_2$ induce a $K_4$. On the other hand $L(x_1)_{E\backslash\{e\}},L(y_1)_{E\backslash\{e\}},L(x_2)_{E\backslash\{e\}},L(y_2)_{E\backslash\{e\}}$ in $\mathcal{M}\backslash\{e\}$ lie together in a rank $2$ contraction minor of $\mathcal{M}\backslash\{e\}$. Thus by Lemma~\ref{lem:rk2}, the vertices $x_1,y_1,x_2,y_2$ should induce a subgraph of a cycle, a contradiction.\\
$(2)\Rightarrow (3):$ Obvious.\\
$(3)\Rightarrow (1):$ We only have to check (C3). So let $L(v)\neq\pm L(w)$ be two labels and $e\in S(L(v),L(w))$. On any $(v,w)$-path $P$ there must be a vertex $u$ with $L(u)_e=0$. If $P$ is crabbed $L(u)$ satisfies (C3) for $L(v),L(w)$ with respect to $e$.
\end{proof}

Cocircuit graphs of general oriented matroids are $2(r-1)$-connected~\cite{Cor-93}. Here we have shown a crabbed analogue for uniform oriented matroids. An interesting question would be if there is a similar result to Theorem~\ref{thm:main} for non-uniform oriented matroids.

Let us now turn to another basic concept for the recognition of cocircuit graphs. A sign labeling $L$ of $G_{\mathcal{M}}$ induces the map $A_L$ which takes $v$ to the unique vertex $w$ with $L(w)=-L(v)$. We call $A_L$ the \textit{AP-labeling} induced by $L$.
\begin{lem}\label{lem:AP}
 If $L$ is an OM-labeling then $A_L$ is an involution in ${\rm{Aut}}(G_{\mathcal{M}})$ which satisfies $d(v,A_L(v))=n-r+2$ for every $v\in V$.
\end{lem}
\begin{proof}
Let $L$ be the OM-labeling inducing the AP-labeling $A_L$. By the definition of sign labeling it is clear that $A_L$ is an automorphism of order $2$. Since any neighbor $u$ of $v$ lies in a coline with $A_L(v)$ and $|S(L(u),-L(v))|+\frac{1}{2}|L^0(u)\Delta L^0(v)|=n-r+1$, by Lemma~\ref{lem:rk2} we have $d(v,A_L(v))=n-r+2$.
\end{proof}

More generally, in a context where the parameters $r,n$ are given an \textit{antipodal labeling} of a graph $G$ is an involution $A\in {\rm{Aut}}(G)$ such that the graph distance satisfies $d(v,A(v))=n-r+2$ for every $v\in V$. 

\section{The Algorithm}

The input is an undirected simple connected graph $G=(V,\mathcal{E})$. The algorithm decides if $G=G_{\mathcal{M}}$ for some {\rm{OM}} $\mathcal{M}=(E,\mathcal{C}^*)$. In the affirmative case it returns $\mathcal{M}$. Otherwise one of the steps of the algorithm
fails.

\begin{enumerate}
 \item Check if $G$ is $2(r-1)$-regular for some $r$.
\item Check if $|V|=2\binom{n}{n-r+1}$ for some $n$. 
\item Calculate the $V\times V$ distance matrix of $G$.
\item Fix $v\in V$ and define $D(v):=\{w\in V\mid d(v,w)=n-r+2\}$.
\item For all $w\in D(v)$ do
\begin{itemize}
	\item[A.] Construct an antipodal labeling $A$ with $A(v)=w$.
	\item[B.] Construct a sign labeling $L$ of $G$ with $A=A_L$.
	\item[C.] Check if $L$ is an OM-labeling. If so, define $\mathcal{C}^*:=L(V)$, return $(E, \mathcal{C}^*)$ and stop.
\end{itemize}
\item Return that $G$ is no cocircuit graph.
\end{enumerate}

Parts 1 and 2 run in time $|\mathcal{E}|=(r-1)|V|$ and are necessary to determine the parameters $r$ and $n$ of $\mathcal{M}$. The distance matrix is computed to avoid repeated application of shortest path algorithms during the main part of the algorithm. Since $G$ is unweighted and undirected we can obtain its distance matrix in $O(|V||\mathcal{E}|)$, see for instance Chapter~6.2 of~\cite{Sch-03}. Hence we can do the
first three parts in $O(r|V|^2)$.

For the rest of the algorithm we have to execute steps A to C at most $|D(v)|$ times. These will be explained in the following.

\subsection{A. Construct an antipodal labeling.}
\begin{lem}
Let $G_{\mathcal{M}}$ be a cocircuit graph with AP-labeling $A_L$. If $A(v)=w$ and $u$ is a neighbor of $v$ then $A(u)$ is the unique neighbor $u'$ of $w$ with $d(u,u')=n-r+2$.
\end{lem}
\begin{proof}
Suppose that, besides $u'=A(u)$, there is another neighbor $u''$ of $w$ with $d(u,u'')=n-r+2$. Since $d(u',u'')\leq 2$, we have $|L^0(u')\Delta L^0(u'')|\leq 4$. Since $L$ is an OM-labeling, $L(u),L(u'),L(u'')$ lie
in a rank $r'=3$ contraction minor on $n'$ elements. Contraction in $\mathcal{M}$ just means deletion of vertices in $G_{\mathcal{M}}$. This implies that we have $d\leq d'$ for the distance functions of the cocircuit graphs of $\mathcal{M}$ and the contraction minor respectively. On the other hand $n'-r'+2=n-(r-3)-3+2=n-r+2$. This yields $d'(u,u')=d'(u,u'')=n-r+2$ in a {\rm{OM}} of rank $3$, which contradicts Lemma~\ref{lem:rk3}.
\end{proof}

We obtain a simple breadth-first search algorithm that given $A(v)=w$ determines $A$ in time $O(r^2|V|)$.
Just walk from the root $v$ through a breadth-first search tree. For vertex $u$ with father $f(u)$ the vertex $A(f(u))$ is known. Look through the $2(r-1)$ neighbors of $A(f(u))$. For the unique neighbor $u'$ having $d(u,u')=n-r+2$ set $A(u)=u'$. If $u'$ is not unique or does not exist, then no AP-labeling $A$ of $G$ with $A(v)=w$ exists.

\subsection{B. Construct a sign labeling of $G$.}
We use the algorithm presented in~\cite{Bab-01}, which given an antipodal labeling $A$ constructs a sign labeling $L$ such that $A=A_L$. If $A$ is an AP-labeling the algorithm finds an OM-labeling witnessing that. Otherwise the algorithm fails. It has runtime $O(rn|V|)$.

\subsection{C. Check if a sing labeling is an OM-labeling.}
We check for every $u\in V(G)$ if there is a crabbed path to every vertex $w\in V(G)$. By Theorem~\ref{thm:main} this is equivalent
to $L$ being an OM-labeling. To improve running time we need the following simple lemma.

\begin{lem}\label{lem:crab}
If a $(u,v)$-path $P$ and a $(v,w)$-path $P'$ are $(L(u),Y)$-crabbed and $(L(v),Y)$-crabbed, respectively, then their concatenation is $(L(u),Y)$-crabbed.
\end{lem}
\begin{proof}
The vertices in $P$ -- in particular $v$ -- satisfy the conditions for being $(L(u),Y)$-crabbed. Hence for every $w\in P'$ we have $L(w)^+ \subseteq L(v)^+ \cup Y^+ \subseteq L(u)^+\cup Y^+$. The analogue statement holds for $L(w)^+$.
\end{proof}

\begin{itemize}
 \item For every $u\in V(G)$ do
\begin{enumerate}
 \item For every edge $\{v,w\}$ do 
\begin{itemize}
\item Delete the undirected edge $\{v,w\}$.
\item If $(v,w)$ is $(L(v),L(u))$-crabbed insert that directed edge.
\item If $(w,v)$ is $(L(w),L(u))$-crabbed insert that directed edge.
\end{itemize}
\item Start a breadth-first search on the resulting directed graph $G'$ at $u$ such that only backward arcs are traversed.
\item If not every vertex is reached by the search, return that $L$ is no OM-labeling and stop.
\end{enumerate}
\item Return that $L$ is an OM-labeling of $G$
\end{itemize}

Lemma~\ref{lem:crab} tells us that for checking if there is a crabbed $(u,v)$-path for every $v\in V(G)$ it is enough to check that if the directed graph $G'$ has a directed path from every vertex to $u$. Step 2 does exactly this. Loop 1 will be executed $(r-1)|V|$ times and each round costs $O(n-r)$ many comparisons. Step 2 runs in time linear in the edges. Since the whole process has to be repeated $|V|$ times, we need $O(r(n-r)|V|^2)$ many operations.

\subsection{Overall Runtime}
We add the runtimes of the single parts of the algorithm. We see that part C dominates all other parts of the algorithm, thus we obtain an overall runtime of $O(|D(v)|r(n-r)|V|^2)$. So far the best known upper bounds for the size of $D(v)$ are in $O(|V|)$, hence our runtime is $O(r(n-r)|V|^3)$. For comparison, the runtime of the algorithm in~\cite{Bab-01} is $O(rn^2|V|^4)$. The improvement of runtime comes from approaching step C in a new way. Already in~\cite{Bab-01} it was asked if in that part some better algorithm was possible.

\section{Antipodality}

The problem of bounding the size of $D(v)$ is hard. In the present section we will point out some open problems concerning this value.
The following lemma is well-known.

\begin{lem}\label{lem:rk3}
If the rank of $\mathcal{M}$ is at most $3$ then for $v,w\in V(G_{\mathcal{M}})$ we have $d(v,w)=n-r+2$ if and only if $-L(v)=L(w)$. Moreover $n-r+2$ is the diameter of $G_{\mathcal{M}}$.
\end{lem}
\begin{proof}
Let $L$ be an OM-labeling of $G_{\mathcal{M}}$ with induced AP-labeling $A_L$. Let $v,w$ be vertices with $A_L(v)\neq w$. We observe the following$\colon$

For any shortest $(v,w)$-path $P$ in $G_{\mathcal{M}}$ we have $P\cap A_L(P)=\emptyset$. On a shortest path there cannot occur anything like $(u',u,\ldots,A_L(u))$, because by Lemma~\ref{lem:rk2}, for every neighbor $u'$ of $u$, we have $d(u',A_L(u))=n-r+1$ since $u'$ and $A_L(u)$ lie in a coline. Hence, $u'$ lies on a shortest $(u,A_L(u))$-path. 

Every shortest path $P=(v=v_0,\ldots,v_k=w)$ satisfies $L(v_i)_e=0$ and $L(v_{i+1})_e\neq 0$, implying $L(w)_e=L(v_{i+1})_e$. Otherwise there would be $v_i,v_j$ in $P$ lying in a coline $\mathcal{M}/e$, but the part of $P$ connecting $v_i,v_j$ would leave $\mathcal{M}/e$. Since $L(v_i)\neq -L(v_j)$, this contradicts Lemma~\ref{lem:rk2}.

This yields that, on a shortest $(v,w)$-path, we will have $L(v_i)_e=0$ and $L(v_{i+1})_e\neq 0$ at most once per $e\notin L^0(v)\cap L^0(w)$ and never if $e\in L^0(v)\cap L^0(w)$. Hence $d(v,w)\leq \underline{L(v)}=n-r+1$.
\end{proof}

The proof actually shows that in an {\rm{OM}} of rank $3$ every shortest path is crabbed. In~\cite{Bab-01} it was asked if the statement of Lemma~\ref{lem:rk3} holds for every rank$\colon$
\begin{quest}\label{quest:1}
Given an {\rm{OM}} of rank $r$ on $n$ elements, does $d(X,Y)=n-r+2$ imply $-X=Y$? Is $n-r+2$ the diameter of $G_{\mathcal{M}}$?
\end{quest}

One could hope that the sign vectors of two vertices $u,v$ give some crucial information about how to connect them by a path. As in the case of rank $3$ matroids one would like to use crabbed paths to prove something about the distance function of $G_{\mathcal{M}}$.
A \textit{tope} in $\mathcal{M}$ is a maximal set $T\subseteq\mathcal{C}^*$ such that $S(X,Y)=\emptyset$ for every $X,Y\in T$. 

\begin{rem}
The assumption that for a fixed tope $T$ of $\mathcal{M}$ every $u,v$ with $L(u),L(v)\in T$ are connected by a crabbed path of length $n-r+1$ implies the Hirsch conjecture.
\end{rem}
\begin{proof}
The Hirsch conjecture says that the graph of a $d$-dimensional simple polytope with $f$ facets has diameter at most $f-d$. Take a $d$-dimensional simple polytope $P$ with $f$ facets in $\mathbb{R}^{d}$. The bounding hyperplanes of $P$ form an affine hyperplane arrangement $H$. Put $H$ into the $(x_{d+1}=1)$-hyperplane of $\mathbb{R}^{d+1}$ and extend $H$ to a central hyperplane arrangement $H'$ in $\mathbb{R}^{d+1}$. We obtain the cone of $P$ as a maximal cell of a central hyperplane arrangement with $f$ hyperplanes. That is, the vertices of $P$ correspond to the cocircuits of a tope $T_P$. This holds for any reorientation of the {\rm{OM}} $\mathcal{M}$ of rank $d$ associated to $H'$. The graph of $P$ is the subgraph of $G_{\mathcal{M}}$ induced by $T_P$. Now for two cocircuits in $T_P$ a path connecting them is crabbed if and only if it is contained in $T_P$.
\end{proof}

The Hirsch conjecture holds in dimension $3$ and topes of {\rm{OM}}s of rank $4$ are combinatorially equivalent to simple polytopes of dimension $3$. Thus, for {\rm{OM}}s of rank $4$ the assumption of the above remark is true for every tope. Still this does not immediately yield a positive answer to Question~\ref{quest:1} for $r=4$. 

Another question, which already seems to be hard and is similar to one posed in~\cite{Bab-01} is the following$\colon$
\begin{quest}\label{quest:2}
How many different antipodal labelings that pass through steps A and B of the algorithm does a cocircuit graph admit?
\end{quest}
Every answer to Question~\ref{quest:2} better than $O(|V|)$ would improve the runtime of our algorithm.

\bibliography{/homes/combi/knauer/Desktop/literature}
\bibliographystyle{amsplain}

\end{document}